\documentclass[12pt]{article}
\usepackage{amsthm,amsmath,amssymb, amsfonts}
\usepackage{enumerate}
\usepackage{url}
\usepackage{hyperref}
\usepackage{fullpage}
\usepackage{fancyhdr}
\usepackage{makeidx}
\usepackage{epic}
\usepackage{eepic}
\usepackage{latexsym}
\usepackage{color}
\usepackage{epsfig}
\usepackage{graphics}
\usepackage{layout}

\newtheorem{theorem}{Theorem}[section]

\newtheorem{lemma}[theorem]{Lemma}
\newtheorem{corollary}[theorem]{Corollary}

\newtheorem{proposition}[theorem]{Proposition}

\newtheorem{example}[theorem]{Example}
\newtheorem{remark}[theorem]{Remark}
\newtheorem*{address}{}
\newtheorem*{email}{}

\numberwithin{equation}{section}

\definecolor{Yellow}{cmyk}{0,0,1,0}
\definecolor{Magenta}{cmyk}{0,1,0,0}
\definecolor{Green}{cmyk}{1,0,1,0.25}
\definecolor{Blue}{cmyk}{1,1,0,0}
\definecolor{Orange}{cmyk}{0,0.5,1,0}
\definecolor{Red}{cmyk}{0,1,1,0}

\renewcommand{\tilde}[1]{\widetilde{#1}}

\hypersetup{colorlinks=true}

\begin {document}

\title{Some explicit computations and models of free products}

\author{Madhushree Basu}


\maketitle


\begin{abstract}
In this note, we first work out some `bare hands' computations of the
most elementary possible free products involving $\mathbb{C}^2 ~(=\mathbb{C} \oplus \mathbb{C} $) and $M_2 ~(= M_2(\mathbb{C}))$. 
Using these, we identify all free products $C \ast D$,
where $C,D$ are of the form $A_1 \oplus A_2$ or $M_2(B)$; $A_1,A_2,B$ are finite von Neumann algebras, as is $A_1 \oplus A_2$
with the `uniform trace'  given by $tr(a_1, a_2) = \frac{1}{2}  (tr(a_1) + tr(a_2))\}$
and $M_2(B)$ with the normalized trace given by $tr((b_{i,j}))=\frac{1}{2}(tr(b_{1,1}) + tr(b_{2,2}))$.
Those results are then used to compute 
various possible free products involving certain finite dimensional von-Neumann algebras,
the free-group von-Neumann algebras and the hyperfinite $II_1$
factor. In the process, we reprove
Dykema's result `$R \ast R \cong LF_2$'. 
\end{abstract}

\section{Introduction}

The motivation behind the computations done in this paper comes from
trying to understand [D1] and [D2].  Most of the results here
have been proved in those two papers.  Further, our proofs seem to
constrain us to only regard direct sums and matrix algebras equipped with the `uniform trace',  and consequently to
results where powers of 2 keep cropping up. The
flip side of the coin is that our proofs involve nothing more than
basic trigonometry and simple cumulant calculations and require little prior background to follow.  

We shall work with the model  $(L^{\infty}([0,\pi/2], \frac{2}{\pi}\int_0^{\frac{\pi}{2}} \cdot dt) \ast
L\mathbb{Z} \ast L\mathbb{Z})$ of $LF_3$. Let $u$ and
$v$ be Haar unitaries generating the two $L\mathbb{Z}$s, and let $c,s
\in L^{\infty}([0,\pi/2])$ be the functions defined by $c(\theta) =
cos ~\theta$ and $s(\theta) = sin ~\theta$ respectively.

Let the trace and cumulant on $M_2(LF_3)$, resp.,
$LF_3$, be denoted by $Tr$ and $\kappa$, resp. $tr$ and $k$ respectively.

Let $U = \begin{pmatrix} 0 & u \\ 0 & 0 \end{pmatrix}$, $V
= \begin{pmatrix} 0 & v \\ 0 & 0 \end{pmatrix}$,  $W
= \begin{pmatrix} c& -s \\ s & c\end{pmatrix}$, and let  
$X = WVW^* = \begin{pmatrix} -cvs & cvc \\ -svs & svc \end{pmatrix}$. 

Then $U$ and $X$ are partial isometries satisfying $U^*U+UU^* = 1$ and
$X^*X+XX^* = 1$; so each of them generates a copy of $M_2$.  

Let $P = UU^* = \begin{pmatrix} 1 & 0 \\ 0 & 0 \end{pmatrix}$ and $Q =
XX^* = \begin{pmatrix} c^2 & cs \\ cs & s^2  \end{pmatrix} = WPW^*$.
Then $P,Q$ are projections of trace $\frac{1}{2}$, each generating a copy of $\mathbb{C}^2$.

Finally, for finite von-Neumann algebra $A,B$ and $k \in \mathbb{N}$, we often use the notations $A^k$ and $M_k(B)$ for
$A \oplus A \oplus \cdots \oplus A (k~ \text{times})$ and
$M_k(\mathbb{C}) \otimes B$ respectively, with the trace always being taken as the `uniform trace'. 

\section{$\mathbb{C}^2\ast\mathbb{C}^2 \cong M_2(L \mathbb{Z})$}

\begin{lemma}\label{c2*c2lem}
\[W^*(\{P, Q \}) = M_2(L \mathbb{Z})\]
\end{lemma}

\begin{proof}
Note that $PQ(1-P) = \begin{pmatrix} 0 & cs \\ 0 &
  0 \end{pmatrix}$. Since $cs$ is positive (in
$L^\infty([0,\frac{\pi}{2}])$) and has no kernel, we see that the
polar decomposition of $PQ(1-P)$ is given by
 $\begin{pmatrix} 0 & cs \\ 0 &
  0 \end{pmatrix} = \begin{pmatrix} 0 & 1 \\ 0 & 0 \end{pmatrix} \begin{pmatrix} 0
  & 0 \\ 0 & cs \end{pmatrix} ~.$ Thus $W^*(\{P, Q \})$ contains all
the matrix units, viz., $P = e_{11}, 1-P = e_{22}, e_{12}=$ polar part of $PQ(1-P)$ and $e_{21}
= e_{12}^*$. On the other hand $W^*(\{P, Q \})$ clearly contains $PQP
= e_{11} \otimes c^2$, and hence also $e_{11} \otimes W^*(\{c^2\})$;
therefore $W^*(\{P, Q \}) \supset e_{11} \otimes
L^\infty([0,\frac{\pi}{2}])$. Finally, by pre- and post-multiplying by
appropriate matrix units, we see that $W^*(\{P, Q \}) \supset e_{ij} \otimes
L^\infty([0,\frac{\pi}{2}]) $ for all $i,j$, and
the proof of the lemma is complete.
\end{proof}

\begin{lemma}\label{PQfree} 
$P$ and $ Q$ are free in $M_2(LF_3)$. 
\end{lemma}

\begin{proof} 
Let $P_0 = P - 1/2 = \begin{pmatrix} 1/2 & 0 \\ 0 &
  -1/2 \end{pmatrix}$ and $Q_0 = Q - 1/2 = \begin{pmatrix} c^2-1/2 &
  cs \\ cs & s^2-1/2 \end{pmatrix}$ be the trace-less versions (i.e.,
translates with trace 0) of $P$ and $Q$. We shall find it convenient
to write $c_n$, resp. $s_n$, for the elements of
$L^\infty([0,\frac{\pi}{2}])$ defined by $c_n(\theta) = cos~n\theta$, resp.,
$s_n(\theta) = sin~n\theta$. 

We are going to verify that the trace of any
alternating product in $2P_0 = \begin{pmatrix} 1 & 0 \\ 0 &
  -1 \end{pmatrix}$ and $2Q_0 =\begin{pmatrix} c_2 & s_2\\
s_2 & -c_2 \end{pmatrix}$is zero. Since we are working with a trace
here, it is enough to prove that $Tr((2P_0.2Q_0)^r) = 0$ and
$Tr((2Q_0.2P_0)^r2Q_0) = 0$. 

However,
\begin{eqnarray*} (2P_0.2Q_0)^r &=& \left( \begin{pmatrix} 1 & 0 \\ 0 &
  -1 \end{pmatrix} \begin{pmatrix} c_2 & s_2 \\ s_2
    & -c_2 \end{pmatrix} \right)^r\\ 
&=& \left( \begin{pmatrix} c_2 & s_2 \\ - s_2
    &  c_2 \end{pmatrix} \right)^r\\ 
&=& \begin{pmatrix} c_{2r} & s_{2r} \\  - s_{2r}
    &  c_{2r} \end{pmatrix}
\end{eqnarray*}

whereas
\begin{eqnarray*} (2Q_0.2P_0)^r2Q_0 &=& \begin{pmatrix} c_{2r} &
    s_{2r} \\  - s_{2r} & c_{2r} \end{pmatrix}  \begin{pmatrix} 1 & 0 \\ 0 &
  -1 \end{pmatrix}\\
&=& \begin{pmatrix} c_{2r+2} & s_{2r+2} \\  s_{2r+2}
    & - c_{2r+2} \end{pmatrix}~,
\end{eqnarray*}

and the desired assertions follows from  

\begin{eqnarray*} 
\frac{2}{\pi}\int_0^{\frac{\pi}{2}}c_{2k}d\theta &=& \frac{2}{\pi}
\frac{1}{2k}(s_{2k}(\pi) - s_{2k}(0))\\
&=& 0 ~ \forall k \in \mathbb{N}
\end{eqnarray*}
\end{proof}

Now the validity of Proposition \ref{c2*c2} follows immediately from
Lemmas \ref{c2*c2lem} and \ref{PQfree}.

\begin{proposition}\label{c2*c2} 
\[\mathbb{C}^2 \ast \mathbb{C}^2 \cong M_2(L \mathbb{Z})\]
\end{proposition}

\section{$M_2 \ast M_2 \cong M_2(LF_3)$}

En route to proving the main result, viz., Proposition \ref{m2*m2}, we
shall also prove Proposition \ref{m2*c2}.

\begin{proposition}\label{m2*c2} 
\[M_2 \ast \mathbb{C}^2 \cong M_2(LF_2)\]
\end{proposition}

\begin{proposition}\label{m2*m2} 
\[M_2 \ast M_2 \cong M_2(LF_3)\]
\end{proposition}

We start with the following crucial lemma.

\begin{lemma}\label{UXfree}
 $U$ and $X$ are $\ast$- free in $M_2(LF_3)$. 
\end{lemma}

\begin{proof}
Observe the following simple facts:

\begin{enumerate}
\item $U^2 = 0 = (U^*)^2, V^2 = 0 = (V^*)^2$ 
\item $P_0^2, Q_0^2 \in \mathbb{C}$ 
\item $PU = U, QX = X, UP = 0, XQ = 0$
\item $P, Q$ are free in $M_2(L^{\infty})$
\end{enumerate}

In view of the the above relations, the sets of trace zero words in
$W^*(\{U\})$ and $W^*(\{X\})$ are linearly generated by $A = \{U, U^*, 2P_0\}$
and $B = \{X, X^*, 2Q_0 \}$ repectively. So we need to check that every
alternating product of elements from $A$ and $B$ has trace zero.  

We may dispose of the `trivial case' when $\Pi$ is an alternating word in
only $P_0$ and $Q_0$, since that is covered by Lemma \ref{PQfree}. 

Consider a typical such product, say $\Pi$; we will prove that
$tr(\Pi_{ij}) = 0$ for all $i,j$ (for the non-trivial case). 
Thus in particular we will have $Tr(\Pi) = \frac{tr(\Pi_{1,1})+tr(\Pi_{2,2})}{2} = 0$ 

One can see that each $\Pi_{ij}$ is a sum of elements
of the form 
\begin{equation}\label{om} 
\omega =  f_0(c,s)w_0f_1(c,s)w_1 f_2(c,s) \cdots w_{n-1}
f_n(c,s) \cdots\end{equation} 
where the $f_i$s are (possibly constant) polynomials  in $\{c,s\}$ and
$w_i \in \{u, u^*, v, v^*\}$. The assumption that we are not dealing
with the trivial case implies that there
must exist at least one $w_i$ in the string in the right hand side of equation \ref{om}.

Now, fix any such word $\omega$ is (of total length $m$, say) as in equation
\ref{om}, which occurs as a summand of $\Pi$; then
\[tr(\omega) = \sum_{\sigma \in NC(m)} k_{\sigma} [f_0(c,s), w_0,
f_1(c,s), \cdots, w_{n-1}, f_n(c,s), \cdots ]~.\] 

We shall show that $k_\sigma (= k_\sigma[f_0(c,s),w_0, f_(c,s),\cdots,w_{n-1},f_n(c,s),
\cdots])= 0$ for each $\sigma \in NC(m)$, for each such $\omega$.

But $u$ and $ v$ are free Haar unitaries and hence R-diagonal. So following 
Proposition 15.1 in \noindent [NS], in order for $k_{\sigma}$ to be possibly non zero, 
it must be the case that
the blocks of $\sigma$ must consist of either only $\{u,u^*\}$ in
alternate positions, or only $\{v, v^*\}$ in alternate positions (in each block 
the number of $u$ = number of $u^*$ and same for $v,v^*$), or
only $\{f_i(c,s):i\}$, all occurring in a non-crossing fashion. Note that in effect
$\omega$ must have the same number of $u$ and $u^*$ as well as the same number of $v$ and $v^*$.

\begin{example}{\rm $\omega_0 = u(cs)v\left(cs(c^2-1/2)s\right)v^*(c^2s)u^* =
  uf_1(c,s)vf_2(c,s)v^*f_3(c,s)u^*$ - with $m=7$ - can possibly give
  non zero cumulant only corresponding to two elements of $ NC(7)$, namely
\[ \left((1,7),(3,5),(2),(4),(6)\right) \mbox{ and }
  \left((1,7),(3,5),(2,6),(4)\right) ~.\] (Here, $f_1(c,s)=cs,
  f_2(c,s)=cs(c^2-1/2)s, f_3(c,s)=c^2s$)}  
\end{example}

As at least one $w_i$ and necessarily also $w_i^*$ must occur in the
string defining $\omega$,  and $\sigma$ is non-crossing,
we see that $k_\sigma$ can hope to be non-zero only if $\omega$ has a
substring of the form $wf(c,s)w^*$, 
with $w \in \{u, u^*, v, v^*\}$ and $f(c,s)$ as above, and with $w$
and $w^*$ in the same block of $\sigma$; thus $\omega$
must contain one of the one of the following four substrings:

\begin{enumerate}
\item $uf(c,s)u^*$
\item $u^*f(c,s)u$
\item $vf(c,s)v^*$
\item $v^*f(c,s)v$
\end{enumerate}

Now consider the above four cases in the following way:
\begin{enumerate}
\item
$uf(c,s)u^*$ can occur as a string in some summand $\omega$
of $\Pi_{ij}$ only if  $USU^*$ occurs as a substring in the alternating
product expression of $\Pi$, where $S =
(2Q_02P_0)^r2Q_0$. Observe in this case that $USU^* = \begin{pmatrix}
  -uc_{2r+2}u^* & 0 \\ 0 & 0 \end{pmatrix}$ and that $f(c,s) = -c_{2r+2}$.
\item Similarly $u^*f(c,s)u$ can occur as a string in some summand $\omega$
of $\Pi_{ij}$ only if  $U^*SU$ occurs as a substring in the alternating
product expression of $\Pi$, where $S =
(2Q_02P_0)^r2Q_0$. Observe in this case that $U^*SU = \begin{pmatrix}
  0 & 0 
  \\ 0 & u^*c_{2r+2}u \end{pmatrix}$ and that $f(c,s) = c_{2r+2}$.
\item Similarly $vf(c,s)v^*$ can occur as a string in some summand $\omega$
of $\Pi_{ij}$ only if  $XSX^*$ occurs as a substring in the alternating
product expression of $\Pi$, where $S =
(2P_02Q_0)^r2P_0$. Observe in this case that $XSX^* = \begin{pmatrix}
  -cvc_{2r+2}v^*c & cvc_{2r+2}v^*s \\ -svc_{2r+2}v^*c &
  -svc_{2r+2}v^*s \end{pmatrix}$ and that again, $f(c,s) = \pm c_{2r+2}$.
\item Similarly $v^*f(c,s)v$ can occur as a string in some summand $\omega$
of $\Pi_{ij}$ only if  $X^*SX$ occurs as a substring in the alternating
product expression of $\Pi$, where $S =
(2P_02Q_0)^r2P_0$. Observe in this case that $X^*SX = \begin{pmatrix}
  sv^*c_{2r+2}vs & -sv^*c_{2r+2}vc \\ -cv^*c_{2r+2}vs &
  cv^*c_{2r+2}vc \end{pmatrix}$ and that again, $f(c,s) = \pm c_{2r+2}$.  
\end{enumerate}

Thus, in any case, we find that for any $\sigma \in NC(m)$  for which
$k_\sigma$ can possibly be non-zero, it must be the case that $\omega$
must contain a substring of the form $w c_{2r+2} w^*$ with $w \in \{u,
u^*, v, v^*\}$ and with $w$ and $w^*$ in the same block of
$\sigma$. Since $W^*(\{c\})$ and $W^*(\{w\})$ are free, we see that the only
way $k_\sigma$ can be non-zero is for $\{c_{2r+2}\}$ to be a block of
$\sigma$; but then $k_\sigma$ has $k(c_{2r+2}) = tr(c_{2r+2}) = 0$ as
a factor. Hence, inded, every $k_\sigma = 0$, as asserted, and the
proof is complete.
\end{proof}

\begin{corollary}\label{PXfree} 
$P$ and $X$ are free in $M_2(LF_3)$. 
\end{corollary}

\begin{proof}
This follows from $P \in W^*(\{U\})$.
\end{proof}

\begin{lemma}\label{m2*c2lem}
\[W^*(\{P,X\}) = M_2(LF_2)\]
\end{lemma}

\begin{proof}
Since $Q = XX^*$, it follows from Lemma \ref{c2*c2lem} that $W^*(\{P,X\}) $ contains each matrix unit $e_{ij}$. 
It follows that $W^*(\{P,X\}) $ contains $M_2(\mathcal{N})$ where $\mathcal{N}$ is the von Neumann algebra generated by the entries of $X$. 
By the last line in the proof of Lemma \ref{c2*c2lem} shows that $L^\infty([0,\frac{\pi}{2}]) \subset \mathcal{N}$. 

Consider the bounded Borel functions 
 $f_n$ defined on $[0,1]$ by $f_n(t) = \left\{ \begin{array}{ll} \frac{1}{t} & t \geq \frac{1}{n}\\0 &  t < \frac{1}{n} \end{array} \right. $ ; 
 observe that $f_n(c)c$ converges strongly to 1 (since $c$ is injective). Since $cvc \in \mathcal{N}$, 
deduce that $v = lim_n (f_n(c)c v cf_n(c)) \in \mathcal{N}$. 
Since $c$ and $v$ are $\ast$-free, and $c,v \in \mathcal{N}$, it follows that $\mathcal{N} \supset LF_2$. 
Since the entries of $P$ and $X$ all lie in $LF_2$, the proof of the Lemma is complete.
\end{proof}

\begin{remark}\label{m2*c2rmk}
 We can prove the above lemma also by proving $W^*(\{U,Q\}) = M_2(LF_2)$ (where $U$ and $Q$ are free),
in an exactly similar way.
\end{remark}

\begin{lemma}\label{m2*m2lem}
\[W^*(\{U,X\}) = M_2(LF_3)\]
\end{lemma}

\begin{proof}
In view of Lemma \ref{m2*c2lem}, we only need to observe that $\{u,v,c\}'' =LF_3$.
\end{proof}

\bigskip
Finally, Proposition \ref{m2*c2} follows from Lemma \ref{m2*c2lem} and Corollary \ref{PXfree}, 
while  Proposition \ref{m2*m2} follows from Lemma \ref{m2*m2lem} and Lemma \ref{UXfree}.

\section{$(A_1 \oplus A_2) \ast (B_1 \oplus B_2) \cong M_2(A_1 \ast A_2 \ast B_1 \ast B_2 \ast L \mathbb{Z})$}

Let $A_1, A_2$ and $B_1, B_2$ be finite von-Neumann algebras. 
We denote the trace and cumulant on matrix algebras over the finite von-Neumann algebras as
 $Tr$ and $\kappa$ and those on the finite von-Neumann algebras themselves as $tr$ and $k$ respectively. 

\bigskip We start by proving a simple but useful lemma on a certain property of the Kreweras compliment
of a non-crossing partition in $NC(n)$ for any $n \in \mathbb{N}$. 

\begin{lemma}\label{krw}
 Let $\pi \in NC(n)$ and $1 \sim_{\pi} n$. Let $V=(k',\cdots,(k+l)')$ be an interval in its Kreweras compliment
$K(\pi)$ for $1 \le k \le n, 0 \le l \le n-k$. 
Then $k \sim_{\pi} (k+l+1)$, where all positive integers are taken modulo $n$.    
\end{lemma}

\begin{proof}
 Note that if $(k)$ is a singloton block then $(k-1)' \sim_{K(\pi)} k'$, a contradiction since $V$ is an interval. 
Similarly $(k+l+1)$ cannot be a singleton block. 
More generally suppose $r_k \in \{1,\cdots,k\}$ and $s_k \in \{k,(k+l+1),\cdots,n\}$ are minimum positive integers such that $k \sim_{\pi} r_k$ and $k \sim_{\pi} s_k$. 
Suppose $r_{k+l+1} \in \{1,\cdots,k,(k+l+1)\}$ and $s_{k+l+1} \in \{(k+l+1),\cdots,n\}$ are the same for $(k+l+1)$.
 We already saw that we cannot have $r_k=s_k=k$ or $r_{k+l+1}=s_{k+l+1}=(k+l+1)$. For simplicity's sake we assume $1 < k, l < (n-k)$.
 The cases $k=1$ or $l=n-k$ will follow similarly.

\bigskip \noindent \textbf{Case $s_k \gneq k$:} In this case $s_k$ being minimum in $\{(k+l+1),\cdots,n\}$ 
such that $k \sim_{\pi} s_k$ we must have $(s_k-1)' \sim_{K(\pi)} (k+l)'$, a contradiction unless $k+l+1 = s_k \sim_{\pi} k$. 

\bigskip \noindent \textbf{Case $s_k = k$:} In this case $r_k$ being minimum in $\{1,\cdots,k-1\}$ 
such that $k \sim_{\pi} r_k$, unless $r_k= 1$ we must have $(r_k-1)' \sim_{K(\pi)} k'$, a contradiction. 
On the other hand if $r_k = 1$, then $r_{k+l+1}$ is forced to be $(k+l+1)$. 
Thus similarly we must have $s_{k+l+1} = n$, otherwise leading into a contradiction.
 But $1 \sim_{\pi} n$. Hence $k \sim_{\pi} r_k = 1 \sim_{\pi} n = s_{k+l+1} \sim_{\pi} (k+l+1)$   
\end{proof}

\begin{proposition}\label{a2*b2}
\[(A_1 \oplus A_2 ) \ast (B_1 \oplus B_2) \cong M_2(A_1 \ast A_2 \ast B_1 \ast B_2 \ast L \mathbb{Z})\]
\end{proposition}

\begin{proof}
Consider the two matrix subalgebras $\begin{pmatrix} A_1 & 0 \\ 0 & A_2 \end{pmatrix} (\cong A_1 \oplus A_2)$ 
and $W \begin{pmatrix} B_1 & 0 \\ 0 & B_2 \end{pmatrix} W^* (\cong B_1 \oplus B_2)$ of $A_1 \ast A_2 \ast B_1 \ast B_2 \ast L \mathbb{Z}$, 
where $W$ is the unitary matrix of the previous sections. 
Following the method used in Lemma \ref{m2*c2lem} we know that these two subalgebras indeed generate the matrix algebra on the right side.
 We need to show that these are free.  

Note that $W \begin{pmatrix} b_1 & 0 \\ 0 & b_2 \end{pmatrix} W^* = \begin{pmatrix} cb_1c + sb_2s & cb_1s-sb_2c \\ sb_1c-cb_2s & sb_1s+cb_2c \end{pmatrix}$

Suppose $\Pi$ is an alternating product of matrices of the form $\begin{pmatrix} a_1 & 0 \\ 0 & a_2 \end{pmatrix}$ and $W \begin{pmatrix} b_1 & 0 \\ 0 & b_2 \end{pmatrix} W^*$, with $a_j \in A_j, b_j \in B_j, tr(a_1+a_2)=0=tr(b_1+b_2)$.

Instead of directly proving $Tr(\Pi)=0$, we will prove a stronger statement:
$tr(\Pi_{i_1,i_1}) = 0 \forall i_1 \in \{1,2\}$. 

We will need to look at alternating products of the form
\begin{align*}
\Pi=&\begin{pmatrix} a_1^1 & 0 \\ 0 & a_2^1 \end{pmatrix} \begin{pmatrix} cb_1^2c + sb_2^2s & cb_1^2s-sb_2^2c \\ sb_1^2c-cb_2^2s & sb_1^2s+cb_2^2c \end{pmatrix} \cdots \begin{pmatrix} a_1^{2r-1} & 0 \\ 0 & a_2^{2r-1} \end{pmatrix} \\
&\begin{pmatrix} cb_1^{2r}c + sb_2^{2r}s & cb_1^{2r}s-sb_2^{2r}c \\ sb_1^{2r}c-cb_2^{2r}s & sb_1^{2r}s+cb_2^{2r}c \end{pmatrix} \begin{pmatrix} a_1^{2r+1} & 0 \\ 0 & a_2^{2r+1} \end{pmatrix},
\end{align*}
where $a_j^i \in A_j, b_j^i \in B_j, tr(a_1^i+a_2^i)=0=tr(b_1^i+b_2^i) \forall i$ - as well as three other kinds of products
 (depending on which sort of matrix the product starts or ends with).

Note that by taking $\begin{pmatrix} a_1^1 & 0 \\ 0 & a_2^1 \end{pmatrix}$ 
or both $\begin{pmatrix} a_1^1 & 0 \\ 0 & a_2^1 \end{pmatrix}$ and $\begin{pmatrix} a_1^{2r+1} & 0 \\ 0 & a_2^{2r+1} \end{pmatrix}$ 
as $\begin{pmatrix}1 & 0 \\ 0 & -1 \end{pmatrix}$, and using the fact that we are working with traces here,
 we find that it is sufficient to consider the one special case listed above. 

For $i_1 \in \{1,2\}$, the $(i_1,i_1)^{th}$ diagonal entry $\Pi_{i_1,i_1}$ of $\Pi$ is a sum of words of the form

$$\omega = a_{i_1}^1 t_{i_1,i_2(\omega)} b_{i_2(\omega)}^2 t'_{i_2(\omega),i_3(\omega)} a_{i_3(\omega)}^3 t_{i_3(\omega),i_4(\omega)} b_{i_4(\omega)}^4 t'_{i_4(\omega),i_5(\omega)} \cdots a_{i_{2r-1}(\omega)}^{2r-1} t_{i_{2r-1}(\omega),i_{2r}(\omega)} b_{i_{2r}(\omega)}^{2r} t'_{i_{2r}(\omega),i_1} a_{i_1}^{2r+1},$$ for $i_2(\omega),\cdots,i_{2r}(\omega) \in \{1,2\}$ and $t_{i_j(\omega),i_{j+1}(\omega)},t'_{i_j(\omega),i_{j+1}(\omega)} \in \{c,s\}$ 
(the reason behind using the cumbersome notation $i_j(\omega)$ is to emphasize the dependence of the indices $i_j$ on the particular summand $\omega$).

Now since the $a$'s and $b$'s are free from the $t$'s and $t^\prime$'s, we see from Theorem 14.4, [NS], 
that $$tr(\omega) = \sum_{\pi \in NC(2r+1)} k_{\pi}(a_{i_1}^1,b_{i_2(\omega)}^2,\cdots,a_{i_1}^{2r+1})tr_{K(\pi)}(t_{i_1,i_2(\omega)},t'_{i_2(\omega),i_3(\omega)},\cdots,t_{i_{2r-1}(\omega),i_{2r}(\omega)},t'_{i_{2r}(\omega),i_1} ,1).$$

Thus, in order to prove that $tr (\Pi_{i_1,i_1})=0$,  it is enough to prove that for any $\pi \in NC(2r+1), i_1 \in \{1,2\}$, 

\begin{equation}\label{etp}
\sum_{\stackrel{\omega \text{ a summand of } \Pi_{i_1,i_1}}{i_2(\omega),\cdots,i_{2r}(\omega) \in \{1,2\}}}k_\pi(a_{i_1}^1,b_{i_2(\omega)}^2,\cdots,a_{i_1}^{2r+1}) tr_{K(\pi)}(t_{i_1,i_2(\omega)},t'_{i_2(\omega),i_3(\omega)}
\cdots,t_{i_{2r-1}(\omega),i_{2r}(\omega)},t'_{i_{2r}(\omega),i_1} ,1) = 0
\end{equation}

The crux of the proof lies in the following key lemma.

\begin{lemma}\label{key}
Let $p \in \mathbb{N}, i_1 \in \{1,2\}$. For $a_i^j \in A_i, b_i^j \in B_k, \mbox{ where }i=1,2, ~j=1,2,\cdots , p+1$, write $a^j = \begin{pmatrix} a_1^{2j-1} & 0 \\ 0 & a_2^{2j-1} \end{pmatrix}$ and $b^j = \begin{pmatrix} cb_1^{2j}c + sb_2^{2j}s & cb_1^{2j}s-sb_2^{2j}c \\ sb_1^{2j}c-cb_2^{2j}s & sb_1^{2j}s+cb_2^{2j}c \end{pmatrix}$; and define 
\[\Pi^\prime = a^1b^1a^2 \cdots a^pb^pa^{p+1}~\mbox{ and }~
\Pi^{\prime\prime} = b^1a^2 \cdots a^{p}b^p.\]

Then
\begin{equation}\label{key1}
\begin{split}
\text{ }\sum_{\stackrel{\omega \text{ a summand of } \Pi'_{i_1,i_1}}{i_2(\omega),\cdots,i_{2p}(\omega) \in \{1,2\}}}&k_{0_{2p+1}}(a_{i_1}^1,b_{i_2(\omega)}^2,\cdots,b_{i_{2p}(\omega)}^{2p},a_{i_1}^{2p+1})tr_{1_{2p+1}}(t_{i_1,i_2(\omega)},t'_{i_2(\omega),i_3(\omega)},\\
&\cdots,t_{i_{2p-1}(\omega),i_{2p}(\omega)},t'_{i_{2p}(\omega),i_1} ,1) = 0
\end{split}
\end{equation}

In particular, when $a^1=a^{p+1} =\begin{pmatrix}1 & 0 \\ 0 & -1 \end{pmatrix}$,
\begin{equation}\label{key2}
\begin{split}
\text{  }\sum_{\stackrel{\omega \text{ a summand of } \Pi''_{i_1,i_1}}{i_2(\omega),\cdots,i_{2p}(\omega) \in \{1,2\}}}&k_{0_{2p-1}}(b_{i_2(\omega)}^2,a_{i_3(\omega)}^3,\cdots,b_{i_{2p}(\omega)}^{2p})tr_{1_{2p}}(t_{i_1,i_2(\omega)},t'_{i_2(\omega),i_3(\omega)},\\
&\cdots,t_{i_{2p-1}(\omega),i_{2p}(\omega)},t'_{i_{2p}(\omega),i_1}) = 0
\end{split}
\end{equation}
\end{lemma}

\begin{proof}
In this proof we write $i_j(\omega)$ as $i_j$ for the sake of simplicity. 
\begin{align*}
&\sum_{\stackrel{\omega \text{ is a summand of } \Pi'_{i_1,i_1}}{i_2,\cdots,i_{2p} \in \{1,2\}}} k_{0_{2p+1}}(a_{i_1}^1,b_{i_2}^2,\cdots,a_{i_1}^{2p+1}) tr_{1_{2p+1}}(t_{i_1,i_2},t'_{i_2,i_3},\cdots,t_{i_{2p-1},i_{2p}},t'_{i_{2p},i_1} ,1)\\
&=\sum_{\stackrel{\omega \text{ is a summand of } \Pi'_{i_1,i_1}}{i_2,\cdots,i_{2p} \in \{1,2\}}} tr(a_{i_1}^1) tr(b_{i_2}^2) \cdots tr(a_{i_1}^{2p+1})tr(t_{i_1,i_2} t'_{i_2,i_3} \cdots t_{i_{2p-1},i_{2p}} t'_{i_{2p},i_1}),
\end{align*}
which is the $(i_1,i_1)$ entry of 

\begin{align*}
&\begin{pmatrix} tr(a_1^1) & 0 \\ 0 & tr(a_2^1) \end{pmatrix} \begin{pmatrix} ctr(b_1^2)c + str(b_2^2)s & ctr(b_1^2)s-str(b_2^2)c \\ str(b_1^2)c-ctr(b_2^2)s & str(b_1^2)s+ctr(b_2^2)c \end{pmatrix} \cdots \begin{pmatrix} tr(a_1^{2p-1}) & 0 \\ 0 & tr(a_2^{2p-1}) \end{pmatrix} \\
&\quad\begin{pmatrix} ctr(b_1^{2p})c + str(b_2^{2p})s & ctr(b_1^{2p})s-str(b_2^{2p})c \\ str(b_1^{2p})c-ctr(b_2^{2p})s & str(b_1^{2p})s+ctr(b_2^{2p})c \end{pmatrix} \begin{pmatrix} tr(a_1^{2p+1}) & 0 \\ 0 & tr(a_2^{2p+1}) \end{pmatrix}=\lambda(2P_02Q_o)^p2P_0 
\end{align*} 

where $\lambda
\begin{cases}
= 0 & \text{ if } \exists i: tr(a_1^i)=tr(a_2^i)=0 \text{ or } tr(b_1^i)=tr(b_2^i)=0 \\
\in \mathbb{C} \setminus \{0\} & \text{ if } \forall i,tr(a_1^i)=-tr(a_2^i) \neq 0, tr(b_1^i)=-tr(b_2^i) \neq 0 
\end{cases}$    

Now the proof follows since by the proof of Lemma \ref{PQfree}, each diagonal entry of an alternating product in $2P_0$ and $2Q_0$ has trace zero.  
\end{proof}

Now let us fix an arbitrary $\pi \in NC(2r+1)$. 

If $\pi=0_{2r+1}$, then equation \ref{etp} follows from equation \ref{key1} for $p=r$.

If $\pi \neq 0_{2r+1}$ then $\exists m < n \in \{1,\cdots,2r+1\}$ such that $m \sim_{\pi} n$. 

Consider an interval, say $V=(k,\cdots, l-1)$ in $K(\pi)$ for $1 \le m \le k \le l-1 \le n \le (r+1)$. Then by Lemma \ref{krw}, $k \sim_{\pi} l$ (since $m \sim_{\pi} n$).

\bigskip \noindent \textbf{Case 1:} Suppose $|V|$ is odd. Then either $a_{i_k}^k$ and $b_{i_l}^l$ or $b_{i_k}^k$ and $a_{i_l}^l$ are joined through $\pi$ 
(depending on whether $k$ or $l$ is odd), which leads to corresponding $k_{\pi}$ being zero, since $A_i,B_j$ are free $\forall i,j$ (Theorem 11.20 [NS]).

\bigskip \noindent \textbf{Case 2:} Suppose $|V|$ is even, both $k$ and $l$ are odd. Then $a_{i_k}^k$ and $a_{i_l}^l$ are joined through $\pi$. Thus $l-1 \gneq k$. 

Also note that if $i_k \neq i_l$ then due to $A_{i_k}$ and $A_{i_l}$ being free the corresponding $k_{\pi}$ would be zero. So we now assume that $i_k = i_l$.

Let $\pi_0 = \pi \setminus \{(k+1) \vee \cdots \vee (l-1)\}$ and $\pi_1 = K(\pi) \setminus \{(k,\cdots, l-1)\}$. 

Then 
\begin{align*}
&\sum_{\stackrel{\omega \text{ is a summand of } \Pi_{i_1,i_1}}{i_2(\omega),\cdots,i_{2r}(\omega)}} k_{\pi}(a_{i_1}^1,b_{i_2(\omega)}^2,\cdots,a_{i_1}^{2r+1}) tr_{K(\pi)}(t_{i_1,i_2(\omega)},\cdots,t'_{i_{2r}(\omega),i_1}) \\
&=\sum_{\stackrel{\omega \text{ is a summand of } \Pi_{i_1,i_1}}{i_2(\omega),\cdots,i_k(\omega)(=i_l(\omega)),i_{l+1}(\omega),\cdots,i_{2r}(\omega)}} k_{\pi_0}(a_{i_1}^1,b_{i_2}^2,\cdots,a_{i_l}^k,a_{i_l}^l,b_{i_{l+1}}^{l+1},\cdots,b_{i_{2r}}^{2r},a_{i_1}^{2r+1})\\
&tr_{\pi_1}(t_{i_1,i_2},\cdots,t'_{i_{k-1},i_k},t_{i_l,i_{l+1}},\cdots,t'_{i_{2r},i_1})\Bigl(\sum_{i_{k+1}(\omega),\cdots,i_{l-1}(\omega)} k_{0_V}(b_{i_{k+1}}^{k+1},\cdots,b_{i_{l-1}}^{l-1})tr_{1_V}(t_{i_k,i_{k+1}},\cdots,t'_{i_{l-1},i_l})\Bigr) \\
&=\sum_{\stackrel{\omega \text{ is a summand of } \Pi_{i_1,i_1}}{i_2(\omega),\cdots,i_k(\omega)(=i_l(\omega)),i_{l+1}(\omega),\cdots,i_{2r}(\omega)}} k_{\pi_0}(a_{i_1}^1,b_{i_2}^2,\cdots,a_{i_l}^k, a_{i_l}^l,b_{i_{l+1}}^{l+1},\cdots,b_{i_{2r}}^{2r},a_{i_1}^{2r+1})\\
&tr_{\pi_1}(t_{i_1,i_2},\cdots,t'_{i_{k-1},i_k},t_{i_l,i_{l+1}},\cdots,t'_{i_{2r},i_1})\Bigl(\sum_{i_{k+1}(\omega),\cdots,i_{l-1}(\omega)} k_{0_{l-k-1}}(b_{i_{k+1}}^{k+1},\cdots,b_{i_{l-1}}^{l-1})tr_{1_{l-k}}(t_{i_k,i_{k+1}},\cdots,t'_{i_{l-1},i_k})\Bigr)
\end{align*}

Write $\Pi=\Pi^{(1,k)}\Pi^{(k+1,l-1)}\Pi^{(l,2r+1)}$, 

where $\Pi^{(1,k)}= a^1b^1\cdots a^\frac{k+1}{2}$, $\Pi^{(k+1,l-1)}= b^\frac{k+1}{2}a^{\frac{k+3}{2}} \cdots b^\frac{l-1}{2}$, $\Pi^{(l,2r+1)}= a^\frac{l+1}{2} \cdots b^{2r}a^{r+1}$.

Further let $\widetilde{\Pi}=\Pi^{(1,k)}\Pi^{(l,2r+1)}$. Note that $\Pi^{(1,k)}$ or $\Pi^{(l,2r+1)}$ can be trivial for the extreme cases.

Then using the fact that $i_k(\omega)=i_l(\omega)$ the above sum may be re-written as

\begin{align*}
&\sum_{\stackrel{\widetilde{\omega} \text{ is a summand of } \widetilde{\Pi}_{i_1,i_1}}{i_2(\widetilde{\omega}),\cdots,i_k(\widetilde{\omega})(=i_l(\widetilde{\omega})),i_{l+1}(\widetilde{\omega}),\cdots,i_{2r}(\widetilde{\omega})}}k_{\pi_0}(a_{i_1}^1,b_{i_2}^2,\cdots,a_{i_l}^k, a_{i_l}^l,b_{i_{l+1}}^{l+1},\cdots,b_{i_{2r}}^{2r},a_{i_1}^{2r+1})\\
&tr_{\pi_1}(t_{i_1,i_2},\cdots,t'_{i_{k-1},i_l},t_{i_l,i_{l+1}},\cdots,t'_{i_{2r},i_1})\Bigl(\sum_{\stackrel{\widetilde{\widetilde{\omega}} \text{ is a summand of } \Pi_{i_k,i_k}^{(k+1,l-1)}}{i_{k+1}(\widetilde{\widetilde{\omega}}),\cdots,i_{l-1}(\widetilde{\widetilde{\omega}})}}k_{0_{l-k-1}}(b_{i_{k+1}}^{k+1},\cdots,b_{i_{l-1}}^{l-1})\\
&tr_{1_{l-k}}(t_{i_k,i_{k+1}},\cdots,t'_{i_{l-1},i_k})\Bigr)
\end{align*}

Now we put $p=\frac{l-k}{2}, \Pi''=\Pi_{i_k,i_k}^{(k+1,l-1)}$ in equation \ref{key2} to get 
for any $i_k(\widetilde{\omega})=i_k(=i_l=i_l(\widetilde{\omega})) \in \{1,2\}$, 
\begin{equation*}
\sum_{\stackrel{\widetilde{\widetilde{\omega}} \text{ is a summand of } \Pi_{i_k,i_k}^{(k+1,l-1)}}{i_{k+1}(\widetilde{\widetilde{\omega}}),\cdots,i_{l-1}(\widetilde{\widetilde{\omega}})}} k_{0_{l-k-1}}(b_{i_{k+1}(\widetilde{\widetilde{\omega}})}^{k+1},\cdots,b_{i_{l-1}(\widetilde{\widetilde{\omega}})}^{l-1})tr_{1_{l-k}}(t_{i_k,i_{k+1}(\widetilde{\widetilde{\omega}})},\cdots,t'_{i_{l-1}(\widetilde{\widetilde{\omega}}),i_k})=0,
\end{equation*}
\noindent
thus proving equation \ref{etp}, as desired. 

\bigskip \noindent \textbf{Case 3:} The case when $|V|$ is even and both $k$ and $l$ are even is proved exactly as in the previous case. 
\end{proof}

\bigskip

\begin{corollary}\label{a2*lz}
\[(A_1 \oplus A_2) \ast L\mathbb{Z} \cong M_2(A_1 \ast A_2 \ast LF_3)\]
\end{corollary}

\begin{proof}
 Follows from Proposition \ref{a2*b2} as well as by noting the fact that $L\mathbb{Z} \cong L\mathbb{Z} \oplus L\mathbb{Z}$, 
both being singly generated von Neumann algebras with non atomic distributions
\end{proof}

\section{$M_2(A) \ast M_2(B) \cong M_2(A \ast B \ast LF_3)$} 

Let $A,B$ be finite von-Neumann algebras. The notations for trace and cumulant remain same as above. 

\begin{proposition}\label{m2a*m2b}
\[M_2(A) \ast M_2(B) \cong M_2(A \ast B \ast LF_3)\]
\end{proposition}

\begin{proof}
Set $Y=\begin{pmatrix}1 & 0 \\ 0 & u\end{pmatrix}$ and $Z=\begin{pmatrix}1 & 0 \\ 0 & v\end{pmatrix}$  in $M_2(LF_3)$, 
where $u,v$ are Haar unitaries as in the introduction, such that $A,B,\{u\},\{v\},\{c\}$ are free. 

We want to show that $Y^*M_2(A)Y (\cong M_2(A))$ and $WZ^*M_2(B)ZW^* (\cong M_2(B))$ are free in $M_2(A \ast B \ast LF_3)$. 
As before, following the method used in Lemma \ref{m2*c2lem}, we can conclude that these two subalgebras generate the algebra on the right in the proposition. 

For $a_{i,j} \in A, b_{i,j} \in B$, 
\begin{align*}
Y^*(a_{i,j})Y &=\begin{pmatrix}1 & 0 \\ 0 & u^* \end{pmatrix} \begin{pmatrix}a_{1,1} & a_{1,2} \\ a_{2,1} & a_{2,2} \end{pmatrix} \begin{pmatrix}1 & 0 \\ 0 & u \end{pmatrix}\\
&= \begin{pmatrix} a_{1,1} & a_{1,2}u \\ u^*a_{2,1} & u^*a_{2,2}u \end{pmatrix}
\end{align*} 
\begin{align*}
WZ^*(b_{i,j})ZW^* &=\begin{pmatrix} c & -s \\ s & c\end{pmatrix} \begin{pmatrix}1 & 0 \\ 0 & v^* \end{pmatrix} \begin{pmatrix}b_{1,1} & b_{1,2} \\ b_{2,1} & b_{2,2} \end{pmatrix}
\begin{pmatrix}1 & 0 \\ 0 & v \end{pmatrix} \begin{pmatrix} c & s \\ -s & c\end{pmatrix} \\
&= \begin{pmatrix} c & -s \\ s & c\end{pmatrix} \begin{pmatrix}b_{1,1} & b_{1,2}v \\ v^*b_{2,1} & v^*b_{2,2}v \end{pmatrix}
\begin{pmatrix} c & s \\ -s & c\end{pmatrix}\\
&= \begin{pmatrix} cb_{1,1}c-sv^*b_{2,1}c-cb_{1,2}vs+sv^*b_{2,2}vs  & cb_{1,1}s-sv^*b_{2,1}s+cb_{1,2}vc-sv^*b_{2,2}vc \\
sb_{1,1}c+cv^*b_{2,1}c-sb_{1,2}vs-cv^*b_{2,2}vs & sb_{1,1}s+cv^*b_{2,1}s+sb_{1,2}vc+cv^*b_{2,2}vc \end{pmatrix}
\end{align*} 

As above, to prove freeness, it is enough to check on an alternating product of the form 
\begin{align*}
\Pi=&\begin{pmatrix} a_{1,1}^1 & a_{1,2}^1u \\ u^*a_{2,1}^1 & u^*a_{2,2}^1u \end{pmatrix} 
\Bigl(\begin{pmatrix} c & -s \\ s & c\end{pmatrix} \begin{pmatrix}b_{1,1}^2 & b_{1,2}^2v\\ v^*b_{2,1}^2 & v^*b_{2,2}^2v \end{pmatrix}
\begin{pmatrix} c & s \\ -s & c\end{pmatrix}\Bigr) \\
&\cdots \Bigl(\begin{pmatrix} c & -s \\ s & c\end{pmatrix} \begin{pmatrix}b_{1,1}^{2r} & b_{1,2}^{2r}v \\ v^*b_{2,1}^{2r} & v^*b_{2,2}^{2r}v \end{pmatrix}
\begin{pmatrix} c & s \\ -s & c\end{pmatrix}\Bigr) \begin{pmatrix} a_{1,1}^{2r+1} & a_{1,2}^{2r+1}u \\ u^*a_{2,1}^{2r+1} & u^*a_{2,2}^{2r+1}u \end{pmatrix}, 
\end{align*}
for $tr(a_{1,1}^i)+tr(a_{2,2}^i)=0=tr(b_{1,1}^j)+tr(b_{2,2}^j)$.

Here too we will prove that each diagonal entry of $\Pi$ has trace zero.

Let $\omega$ be a summand of $(i_1,i_1)^{th}$ diagonal entry of the above product. Then $\omega$ is an alternating product of the form
\begin{equation*}
a_{i_1,i_2}^1w_{i_2,i_3}^2b_{i_3,i_4}^3w_{i_4,i_5}^4 \cdots w_{i_{2r-2},i_{2r-1}}^{2r-2}b_{i_{2r-1},i_{2r}}^{2r-1}w_{i_{2r},i_{2r+1}}^{2r}a_{i_{2r+1},i_1}^{2r+1} 
\end{equation*}
where $w_{i_k,i_{k+1}}^k \in \pm \{c,sv^*,us, ucv^*,vs,su^*,vcu^*\}$ depending on $\omega$.

As before from [NS] we can say that $$tr(\omega)=\sum_{\pi \in NC(r+1)}k_{\pi}(a_{i_1,i_2}^1,b_{i_3,i_4}^3,\cdots,
b_{i_{2r-1},i_{2r}}^{2r-1},a_{i_{2r+1},i_1}^{2r+1})tr_{K(\pi)}(w_{i_2,i_3}^2,w_{i_4,i_5}^4,\cdots,w_{i_{2r},i_{2r+1}}^{2r},1).$$

Using the fact that $u^*u=v^*v=1$, as in Lemma \ref{key}, here also
$$\sum_{\omega \text{ is a summand of }
\Pi_{i_1,i_1}} k_{0_{r+1}}(a_{i_1,i_2}^1,b_{i_3,i_4}^3, \cdots, b_{i_{2r-1},i_{2r}}^{2r-1},a_{i_{2r+1},i_1}^{2r+1})tr_{1_{r+1}}(w_{i_2,i_3}^2,w_{i_4,i_5}^4,\cdots,w_{i_{2r},i_{2r+1}}^{2r},1)$$
is the $(i_1,i_1)^{th}$ entry of the matrix 

\begin{align*}
\Bigl(tr(a_{1,2}^1)U&+tr(a_{2,1}^1)U^*+tr(a_{1,1}^1) 2P_0\Bigr) \Bigl(W(tr(b_{1,2}^2)V+tr(b_{2,1}^2)V^*+tr(b_{1,1}^2) 2P_0)W^*\Bigr) \\
&\cdots \Bigl(W(tr(b_{1,2}^{2r})V+tr(b_{2,1}^{2r})V^*+tr(b_{1,1}^{2r}) 2P_0)W^*\Bigr)\Bigl(tr(a_{1,2}^{2r+1})U+tr(a_{2,1}^{2r+1})U^*+tr(a_{1,1}^{2r+1})2P_0\Bigr)\\
=\Bigl(tr(a_{1,2}^1)U&+tr(a_{2,1}^1)U^*+tr(a_{1,1}^1)2P_0\Bigr)\Bigl(tr(b_{1,2}^2)X+tr(b_{2,1}^2)X^*+tr(b_{1,1}^2)2Q_0\Bigr)\\
&\cdots \Bigl(tr(b_{1,2}^{2r})X+tr(b_{2,1}^{2r})X^*+tr(b_{1,1}^{2r})2Q_0\Bigr)\Bigl(tr(a_{1,2}^{2r+1})U+tr(a_{2,1}^{2r+1})U^*+tr(a_{1,1}^{2r+1})2P_0\Bigr),
\end{align*} 

where $U,V,X$ are trace zero partial isometries in $M_2(LF_3)$ as defined in the introduction,

But by proof of Lemma \ref{UXfree}, each diagonal entry of an alternating product in $\{U,U^*,2P_0\}$ and $\{X,X^*,2Q_0\}$ has trace zero.

Thus exactly as Lemma \ref{key}, for any $i_1 \in \{1,2\}$, $$\sum_{\omega \text{ is a summand of }\Pi_{i_1,i_1}} k_{0_{r+1}}(a_{i_1,i_2}^1,b_{i_3,i_4}^3, \cdots, b_{i_{2r-1},i_{2r}}^{2r-1},a_{i_{2r+1},i_1}^{2r+1})tr_{1_{r+1}}(w_{i_2,i_3}^2,w_{i_4,i_5}^4,\cdots,w_{i_{2r},i_{2r+1}}^{2r},1)=0.$$ 

Now rest of the proof follows similarly as Proposition \ref{a2*b2}.  
\end{proof}


\bigskip

\begin{proposition}\label{a2*m2b} 
\[(A_1 \oplus A_2) \ast M_2(B) \cong M_2(A_1 \ast A_2 \ast B \ast LF_2)\]
\end{proposition}

\begin{proof}
Here one needs to prove that $\begin{pmatrix}A_1 & 0 \\ 0 & A_2 \end{pmatrix} (\cong A_1 \oplus A_2)$ and $WZ^*M_2(B)ZW^* (\cong M_2(B))$
 are free and they generate RHS. The proof is exactly similar to that of Proposition \ref{a2*b2} or Proposition \ref{m2a*m2b}, using Corollary \ref{PXfree}. 
\end{proof}

\bigskip

\begin{remark}\label{a2*m2brmk}
 As in Remark \ref{m2*c2rmk} here too we can have an exactly similar alternate proof using $W \begin{pmatrix} A_1 & 0 \\ 0 & A_2 \end{pmatrix} W^*$
and $Y M_2(B) Y^*$ as model.                                                                                                                                     
\end{remark}

\bigskip

\begin{corollary}\label{m2a*lz}
\[M_2(A) \ast L\mathbb{Z} \cong M_2(A \ast LF_4)\] 
\end{corollary}

\bigskip

\begin{corollary}\label{cor}
For $k_1,k_2,l_1,l_2 \in \mathbb{N} \cup \{0\}$,
\begin{enumerate}
\item $(LF_{k_1} \oplus LF_{k_2}) \ast (LF_{l_1} \oplus LF_{l_2} ) \cong M_2(LF_{k_1+k_2+l_1+l_2+1})$
\item $(LF_{k_1} \oplus LF_{k_2}) \ast L\mathbb{Z} \cong M_2(LF_{k_1+k_2+3})$
\item $M_2(LF_k) \ast M_2(LF_l) \cong M_2(LF_{k+l+3})$
\item $(LF_{k_1} \oplus LF_{k_2}) \ast M_2(LF_l) \cong M_2(LF_{k_1+k_2+l+2})$
\item $L\mathbb{Z} \ast M_2(LF_l) \cong M_2(LF_{l+4})$
\end{enumerate}
\end{corollary}

\begin{proof}
(1), (3) and (4) are direct consequences of Proposition \ref{a2*b2}, Proposition \ref{m2a*m2b} and Proposition \ref{a2*m2b} respectively. 
(2) and (5) follow from those as well as Corollary \ref{a2*lz} and \ref{m2a*lz}
(In fact (3) follows directly from Theorem 5.4.1 [VDN]).
\end{proof}


\section{Applications}

In this section, we use the results proved in the previous sections and compute various possible free products involving
the hyperfinite $II_1$ factor $R$,
$(LF_l)^{2^m}$ and $M_{2^n}(LF_k)$, where $m,n \in \mathbb{N}, k,l \in \mathbb{N} \cup \{0\}$ ($LF_0$ is considered as $\mathbb{C}$). 
These results were proved in section 1-3 of [D2] in a much more general context, but with a different approach.

For the computation we frequently need to use
Theorem 5.4.1 of [VDN]. For our purpose it will be enough to stick to the following 2-dimensional version of the theorem: For $k \ge 2$,

\begin{align}\label{voi}
 LF_k \cong M_2(LF_{(k-1)4+1}).
\end{align}


\bigskip
 
We now state an example of computing a free product of certain finite dimensional algebras, using the above sections:


\begin{example}\label{c2n*c2m}
$\mathbb{C}^{2^n} \ast \mathbb{C}^{2^m} \cong M_2(LF_{5-2(\frac{1}{2^{n-1}}+\frac{1}{2^{m-1}})}), n \ge 1$. In particular,
$\mathbb{C}^{2^n} \ast \mathbb{C}^{2^n} \cong M_2(LF_{5-\frac{4}{2^{n-1}}}), n \ge 1$.
\end{example}

\begin{proof}
We first prove inductively  that $\mathbb{C}^{2^n} \ast \mathbb{C}^{2^n} \cong M_2(LF_{a_n})$ for some $a_n \in [1, \infty)$, for all $n \ge 1$. Then the basic step follows from Proposition \ref{c2*c2}. Also $a_1=1$ by the same Proposition.

\begin{align*}
\mathbb{C}^{2^{n+1}} \ast \mathbb{C}^{2^{n+1}} &\cong (\mathbb{C}^{2^n} \oplus \mathbb{C}^{2^n}) \ast (\mathbb{C}^{2^n} \oplus \mathbb{C}^{2^n}) \\
&\cong M_2(\mathbb{C}^{2^n} \ast \mathbb{C}^{2^n} \ast \mathbb{C}^{2^n} \ast \mathbb{C}^{2^n} \ast L \mathbb{Z}), \text{by Proposition \ref{c2*c2}} \\
&\cong M_2(M_2(LF_{a_n}) \ast M_2(LF_{a_n}) \ast L\mathbb{Z}), \text{by induction hypothesis}  \\
&\cong M_2(M_2(LF_{2a_n+3}) \ast L\mathbb{Z}), \text{by Corollary \ref{cor}}  \\
&\cong M_2(LF_{\frac{2a_n+3-1}{4}+2}\ast L\mathbb{Z}), \text{by equation \ref{voi}}  \\
&\cong M_2(LF_{\frac{a_n}{2}+\frac{5}{2}})
\end{align*}

Thus the induction is complete. Moreover we have the recurrence relation $a_{n+1} = \frac{a_n}{2}+\frac{5}{2}$.

Now,
\begin{align*}
a_{n+1} &= \frac{a_n}{2} + \frac{5}{2} = \frac{\frac{a_{n-1}}{2}+\frac{5}{2}}{2}+\frac{5}{2}\\
&= \frac{a_{n-1}}{2^2}+\frac{5}{2.2}+\frac{5}{2} = \frac{a_1}{2^n}+\frac{5}{2}(\frac{1}{2^{n-1}}+\frac{1}{2^{n-2}}+\cdots+1) \\
&= \frac{1}{2^n}+\frac{5}{2}(\frac{1}{2^{n-1}}+\frac{1}{2^{n-2}}+\cdots+1) = \frac{1}{2^n}+\frac{5(2^n-1)}{2^n} = 5-\frac{4}{2^n}.
\end{align*}

From the above calculations and equation \ref{voi}, 

\begin{align*}
 \mathbb{C}^{2^n} \ast \mathbb{C}^{2^n} &\cong M_2(LF_{5-\frac{4}{2^{n-1}}}),\text{ as required} \\
&\cong LF_{2-\frac{1}{2^{n-1}}}, n \ge 1.
\end{align*}

Without loss of generality we can assume that $n \ge m$. For $m=n=1$ the proof follows from Proposition \ref{c2*c2}.

Let $n \ge 2$. Then,
 \begin{align*}
  \mathbb{C}^{2^n} \ast \mathbb{C}^{2^m} &\cong M_2(\mathbb{C}^{2^{n-1}} \ast \mathbb{C}^{2^{n-1}} \ast \mathbb{C}^{2^{m-1}}  \ast \mathbb{C}^{2^{m-1}} \ast L\mathbb{Z}) \\
  &\cong M_2(LF_{2-\frac{1}{2^{n-2}}} \ast LF_{2-\frac{1}{2^{m-2}}} \ast L\mathbb{Z})\\
  &\cong M_2(LF_{5-2(\frac{1}{2^{n-1}}+\frac{1}{2^{m-1}})}).
 \end{align*}
\end{proof}

We are now ready to state the following proposition that summarizes the promised computations of the free products involving certain finite dimensional
von-Neumann algebras and the free-group von-Neumann algebras.
We will omit the proof since it is a simple exercise of induction using the previous sections,
similar to the above example.

\vspace{7 mm}

\begin{proposition}\label{thm}
For $m,n \in \mathbb{N}$, $k,l \in \mathbb{N} \cup \{0\}$ and $LF_0=\mathbb{C}$, 
\begin{enumerate}
\item $(LF_k)^{2^n} \ast (LF_l)^{2^m} \cong M_2(LF_{5+\frac{2(k-1)}{2^{n-1}}+\frac{2(l-1)}{2^{m-1}}})$. 
In particular $(L\mathbb{Z})^{2^n } \ast (L\mathbb{Z})^{2^m} \cong M_2(LF_5)$.
\item $M_{2^n}(LF_k) \ast (LF_l)^{2^m} \cong$ 
$M_2(LF_{5+\frac{k-1}{4^{n-1}}+\frac{2(l-1)}{2^{m-1}}})$.
In particular $M_{2^n}(L\mathbb{Z}) \ast (L\mathbb{Z})^{2^m} \cong M_2(LF_5)$.
\item $M_{2^n}(LF_k) \ast M_{2^m}(LF_l) \cong M_2(LF_{5+\frac{k-1}{4^{n-1}}+\frac{l-1}{4^{m-1}}})$.
In particular $M_{2^n}(L\mathbb{Z}) \ast M_{2^m}(L\mathbb{Z}) \cong M_2(LF_5)$. 
\end{enumerate}
\end{proposition}

\bigskip

\begin{remark} 
There exist explicitly computable functions $f,g: [1, \infty) \rightarrow [1, \infty)$,
such that  whenever finite von-Neumann algebras $A, B$ satisfy $A \ast B \cong LF_t$ for some $t \in [1, \infty)$, then, for all $n \in \mathbb{N}$, we have 

\begin{itemize}
\item $A^{2^n} \ast B^{2^n} \cong M_2(LF_{f(t)})$, 
\item $M_{2^n}(A) \ast M_{2^n}(B) \cong M_2(LF_{g(t)})$. 
\end{itemize}
\end{remark}

\begin{remark} 
 One can obviously extend the above proposition by taking $LF_{k_1} \oplus \cdots \oplus LF_{k_{2^n}}$ for $k_i \ge 0$, instead of $(LF_k)^{2^n}$.
\end{remark}

\bigskip

We know that the \textbf{hyperfinite $II_1$ factor $R$} can be constructed as an infinite tensor product of type $I_{2^n}$ factors,
i.e. scalar matrix algebras of dimension $2^n$. Again $L\mathbb{Z} \cong L^{\infty}([0,\frac{\pi}{2}])$
can be thought of as an infinite tensor product of $\mathbb{C}^{2^n}$.

Proposition (\ref{thm}) suggests - on `taking the limit as $m,n \rightarrow \infty$' , with $k=l=0$ - that
\begin{enumerate}
\item $L\mathbb{Z} \ast L\mathbb{Z} \cong M_2(LF_5)$ (trivially true);
\item $R \ast L\mathbb{Z} \cong M_2(LF_5)$ (see Theorem 5.4.3 [VDN]);
\item $R \ast R \cong M_2(LF_5)$.
\end{enumerate}

We could not come up with a matrix model to prove the above statements, approximating $L\mathbb{Z}$ and $R$ as
by finite dimensional algebras as $n \rightarrow \infty$ in Proposition \ref{thm}.
But we shall indeed give a rigorous proof for the assertion about $R \ast R$, as against the `limiting' heuristics:
 
In view of the uniqueness of the hyperfinite $II_1$ factor (see [MvN]), we know that

\begin{align}\label{rm2r}
R \cong M_2(R)
\end{align}

Now using Theorem 5.4.3 of [VDN], i.e.
\begin{align*}\label{r*lfk}
LF_k \ast R \cong LF_{k+1}
\end{align*} 
we may deduce the following:

\vspace{5 mm}

\begin{proposition}\label{r}
For finite von-Neumann algebra
$A_1,A_2,B$,
\begin{enumerate}
\item $R \ast (A_1 \oplus A_2) \cong M_2(A_1 \ast A_2 \ast LF_3)$
\item $R \ast M_2(B) \cong M_2(B \ast LF_4)$
\end{enumerate}\
\end{proposition}

\begin{proof}
By above equations and Proposition \ref{a2*m2b}
\begin{align*}
R \ast (A_1 \oplus A_2) &\cong M_2(R) \ast (A_1 \oplus A_2) \\
&\cong M_2(A_1 \ast A_2 \ast R \ast L\mathbb{Z}) \\
&\cong M_2(A_1 \ast A_2 \ast LF_3) 
\end{align*}

The other statement follows similarly using Proposition \ref{m2a*m2b}.
\end{proof}

\bigskip	

\begin{corollary}\label{r*r}
$R \ast R \cong M_2(LF_5)$
\end{corollary}

\begin{proof}
It follows from equation \ref{rm2r} and the above proposition.
\end{proof}

We now wish to observe that our proof can also be led to the strengthened version
Proposition \ref{r*rdyk} of Corollary \ref{r*r}). Following [D1], let us write the left side as
$R \ast \tilde{R}$ to emphasize the distinction between the two free copies of the hyperfinite $II_1$ factor.
Consider $R \cong M_2(R), \tilde{R} \cong M_2(\tilde{R})$.
Then we notice that by proof of Proposition \ref{m2a*m2b}, $M_2(R)$ on the left hand side
gets mapped into $M_2(R)$ on the right hand side (which is $M_2(R \ast \tilde{R} \ast LF_3) \cong M_2(R \ast LF_4)$)
as conjugated by the unitary matrix
$Y^*$ where $Y=\begin{pmatrix} 1 & 0 \\ 0 & u \end{pmatrix}$.

On the other hand, by Proposition \ref{a2*m2b} and Theorem 5.4.3 of [VDN], we have
\begin{align*}
 M_2(R) \ast L\mathbb{Z} &\cong M_2(R) \ast (L\mathbb{Z} \oplus L\mathbb{Z})\\
&\cong M_2(R \ast LF_4),
\end{align*}  
where by Remark \ref{a2*m2brmk}, $M_2(R)$ on the left hand side is mapped into $M_2(R)$ on the right hand side
in exactly the same manner as above, i.e. as conjugated by the same unitary matrix $Y^*$.

In fact we note here that for projection $P=\begin{pmatrix}1 & 0 \\ 0 & 0 \end{pmatrix} \in M_2(R)$,

\begin{align*}
P(M_2(R) \ast M_2(\tilde{R}))P \cong PM_2(R)P \ast LF_4 \cong P(M_2(R) \ast L\mathbb{Z})P,
\end{align*}  
where the isomorphisms restricted to $PM_2(R)P$ (which is naturally isomorphic to $R$), in all three cases are the identity maps.

Thus similarly as in Corollary 3.6 of [D1], we can conclude that

\begin{proposition}\label{r*rdyk}
 \[R \ast \tilde{R} \cong R \ast L\mathbb{Z} \cong M_2(LF_5) \cong LF_2,\]
where the first isomorphism restricted to $R$ on the left hand side is the identity map to $R$ on the right hand side. 
\end{proposition}

\section{Acknowledgments} 

I wish to thank V. S. Sunder and Vijay Kodiyalam for helpful discussions. 

\section{References}

\noindent [D1] K. Dykema: Interpolated  factor Factors, Pacific J. Math. Volume 163, Number 1 (1994), 123-135

\noindent [D2] K. Dykema: Free products of hyperfinite von Neumann algebras and free dimension, Duke Math. J. Volume 69, Number 1 (1993), 97-119. 

\noindent [NS] A. Nica, R. Speicher: Lectures on the Combinatorics of Free Probability. Cambridge University Press, 2006

\noindent [VDN] D. V. Voiculescu, K. Dykema, A. Nica: Free Random Variables, CRM Monograph Series, Volume 1, American Mathematical Society, Providence, RI, 1992.

\noindent [MvN] F.J. Murray, J. von Neumann: On rings of operators IV Ann. of Math. (2) , 44 (1943) pp. 716808.

\bigskip


\address{\small INSTITUTE OF MATHEMATICAL SCIENCES, CIT Campus, Taramani, Chennai 600113, India.}                                                                                                   

\email{\small E-mail address: madhushree@imsc.res.in}

\end{document}